\documentclass[a4paper,12pt]{article}
\usepackage[latin1]{inputenc}
\usepackage{amssymb}
\title{{\bf  An Example of  ${\bf \Pi^0_3}$-complete
\\ Infinitary Rational Relation}
 \footnote{LIP Research Report RR 2007-09} }
\author{Olivier Finkel\\{\it Equipe Modèles de Calcul et Complexité}  
 \\ {\it Laboratoire de l'Informatique du Parallélisme}\footnote{UMR 5668 - CNRS - ENS Lyon - UCB Lyon - INRIA }
 \\  CNRS et Ecole Normale Supérieure de Lyon
 \\ 46, Allée d'Italie 69364 Lyon Cedex 07, France.\\ Olivier.Finkel@ens-lyon.fr }
\date{}
\begin{document}

\newtheorem{The}{Theorem}[section]
\newtheorem{Pro}[The]{Proposition}
\newtheorem{Deff}[The]{Definition}
\newtheorem{Lem}[The]{Lemma}
\newtheorem{Rem}[The]{Remark}
\newtheorem{Exa}[The]{Example}
\newtheorem{Cor}[The]{Corollary}

\newcommand{\fa}{\forall}
\newcommand{\Ga}{\Gamma}
\newcommand{\Gas}{\Gamma^\star}
\newcommand{\Gao}{\Gamma^\omega}

\newcommand{\Si}{\Sigma}
\newcommand{\Sis}{\Sigma^\star}
\newcommand{\Sio}{\Sigma^\omega}
\newcommand{\ra}{\rightarrow}
\newcommand{\hs}{\hspace{12mm}

\noi}
\newcommand{\lra}{\leftrightarrow}
\newcommand{\la}{language}
\newcommand{\ite}{\item}
\newcommand{\Lp}{L(\varphi)}
\newcommand{\abs}{\{a, b\}^\star}
\newcommand{\abcs}{\{a, b, c \}^\star}
\newcommand{\ol}{ $\omega$-language}
\newcommand{\orl}{ $\omega$-regular language}
\newcommand{\om}{\omega}
\newcommand{\nl}{\newline}
\newcommand{\noi}{\noindent}
\newcommand{\tla}{\twoheadleftarrow}
\newcommand{\de}{deterministic }
\newcommand{\proo}{\noi {\bf Proof.} }
\newcommand {\ep}{\hfill $\square$}

\maketitle

\begin{abstract}
\noi  We give in this paper an example of  infinitary rational relation,  accepted 
by a 2-tape B\"uchi automaton,  which is ${\bf \Pi^0_3}$-complete in the Borel hierarchy. 
Moreover the example of infinitary rational relation 
given in this paper has a very simple structure and can be easily 
described by its sections.  

\end{abstract}

\noi {\small {\bf Keywords:} infinitary rational relations; topological properties; Borel hierarchy; ${\bf \Pi^0_3}$-complete set.}

\section{Introduction}

Acceptance of infinite words by finite automata was firstly 
considered in the sixties by B\"uchi in order to study 
decidability of the monadic second order theory 
of one successor over the integers \cite{bu62}. Then the so called \orl s have been 
intensively studied and many applications have been found, see \cite{tho,sta,pp} for many results and references. 
\nl Since then many extensions of \orl s have been investigated as the classes of 
\ol s accepted by pushdown automata, Petri nets, Turing machines, 
see \cite{tho,eh,sta} for a survey of this work. 
\nl On the other side
rational relations on finite words were studied in the sixties and played 
a fundamental role in the study of families of context free languages \cite{ber}. 
Investigations on their extension to rational  relations on infinite words were carried out 
or mentioned in the books \cite{bt,ls}. Gire  and Nivat 
studied infinitary rational relations in  \cite{gire1,gn}. 
Infinitary rational relations 
are subsets of  $\Si_1^\om \times \Si_2^\om \times \ldots \times \Si_n^\om$,  
where $n$ is an integer $\geq 2$ and  
$\Si_1$, $\Si_2$, \ldots $\Si_n$   are finite alphabets, which are accepted  by 
 $n$-tape finite B\"uchi automata with $n$  asynchronous 
reading heads. 
So the class of infinitary rational relations extends both the 
class of finitary rational relations {\bf and } the class of \orl s.  
\nl They have been much studied, in particular in connection with 
the rational functions they may define, see for example  \cite{cg,bcps,sim,sta,pri} for many results and references. 
\nl Notice that 
a rational relation $R\subseteq \Si_1^\om \times \Si_2^\om \times \ldots \times \Si_n^\om$ 
may be seen as 
an \ol~ over the product alphabet $\Si_1 \times \Si_2 \times \ldots \times \Si_n$. 
\nl  A way to study the complexity of languages of infinite words  
accepted by finite machines is to study their topological complexity and firstly
to locate them with regard to 
the Borel and the projective hierarchies. 
This work is analysed 
for example in \cite{stac,tho,eh,lt,sta}.  
It is well known that every \ol~ accepted by a Turing machine with a 
B\"uchi or Muller acceptance condition is an analytic set and 
 that \orl s are boolean combinations of ${\bf \Pi_2^0}$-sets 
hence ${\bf \Delta_3^0}$-sets,  \cite{sta,pp}.  
\nl  The question of the topological complexity of  relations on infinite words also 
naturally arises and was asked by Simonnet in \cite{sim}. It was also posed in a more 
general form by Lescow and Thomas in \cite{lt} 
(for infinite labelled partial orders) and in \cite{tho89} 
where Thomas suggested to study reducibility notions and associated completeness results.  
\nl Every infinitary rational relation is an analytic set. 
 We showed in \cite{relrat} that there exist some infinitary rational relations 
which are analytic but non  Borel sets. 
Considering Borel infinitary rational relations we prove in this paper 
that there exist some infinitary rational relations, accepted 
by 2-tape B\"uchi automata, 
which are  ${\bf \Pi^0_3}$-complete. 
\nl  Examples of ${\bf \Si^0_3}$-complete and ${\bf \Pi^0_3}$-complete 
infinitary rational relations have already been given in the conference paper 
\cite{relratbor}. But the proof of the existence of 
${\bf \Pi^0_3}$-complete infinitary rational relations was only sketched and 
we used a coding of 
$\om^2$-words by pairs of infinite words. 
We use in this paper a different coding of $\om^2$-words. This way we get some 
infinitary rational relations which have  a very simple structure and 
can be easily described by their sections.  
\nl The result given in this paper has two interests: 
1) It gives a complete proof of a result of \cite{relratbor}. 
2) Some new ideas are here introduced with a new coding of $\om^2$-words. Some of these new ideas led us 
further to the proof of very surprising results, answering the  long standing open questions of the topological complexity 
of context free $\om$-languages and of infinitary rational relations. In particular infinitary rational relations have the same 
topological complexity as $\om$-languages accepted by B\"uchi Turing machines \cite{cie05,stacs06}and for every recursive ordinal 
$\alpha$ there exist some ${\bf \Pi^0_\alpha}$-complete and some ${\bf \Si^0_\alpha}$-complete infinitary rational relations. 
\nl The result presented in this paper is still interesting, although the result of the paper \cite{stacs06} is stronger; we use  
here a  coding of $\om^2$-words while in \cite{cie05,stacs06} we used a simulation of Turing machines and the examples of infinitary 
rational relations we obtain  are different. 
\nl The result of this paper may also be compared with examples of ${\bf \Si^0_3}$-complete \ol s accepted 
by deterministic pushdown automata with the acceptance condition: ``some stack content 
 appears infinitely often during an infinite run",   
given by Cachat, Duparc,  and Thomas in \cite{cdt} 
or with examples 
of ${\bf \Si^0_n}$-complete and ${\bf \Pi^0_n}$-complete  \ol s,  $n\geq 1$,  accepted  by  
non-deterministic pushdown automata with B\"uchi acceptance condition given in \cite{fina}. 
\nl  The paper is organized as follows. In section 2 we introduce 
the notion of infinitary rational relations. In section 3 we recall definitions of Borel
sets, and we prove our main result  in section 4.

\section{Infinitary rational relations}

\noi  Let $\Si$ be a finite alphabet whose elements are called letters.
A non-empty finite word over $\Si$ is a finite sequence of letters:
 $x=a_1a_2\ldots a_n$ where for all integers $i\in [1; n]$ $a_i \in\Si$.
 We shall denote $x(i)=a_i$ the $i^{th}$ letter of $x$
and $x[i]=x(1)\ldots x(i)$ for $i\leq n$. The length of $x$ is $|x|=n$.
The empty word will be denoted by $\lambda$ and has 0 letter. Its length is 0.
 The set of finite words over $\Si$ is denoted $\Sis$.
 A (finitary) language $L$ over $\Si$ is a subset of $\Sis$.
 The usual concatenation product of $u$ and $v$ will be denoted by $u.v$ or just  $uv$.
 For $V\subseteq \Sis$, we denote  \quad 
$V^\star = \{ v_1\ldots v_n  \mid ~~\fa i \in [1; n] ~~~ 
   v_i \in V  ~~ \} \cup \{\lambda\}$.

\hs The first infinite ordinal is $\om$.
An $\om$-word over $\Si$ is an $\om$ -sequence $a_1a_2\ldots a_n \ldots$, where 
for all integers  $i\geq 1$~~ $a_i \in\Sigma$.  
 When $\sigma$ is an $\om$-word over $\Si$, we write
 $\sigma =\sigma(1)\sigma(2)\ldots  \sigma(n) \ldots $
and $\sigma[n]=\sigma(1)\sigma(2)\ldots  \sigma(n)$ the finite word of length $n$, 
prefix of $\sigma$.
The set of $\om$-words over  the alphabet $\Si$ is denoted by $\Si^\om$.
 An  $\om$-language over an alphabet $\Sigma$ is a subset of  $\Si^\om$.
For $V\subseteq \Sis$,
 $V^\om = \{ \sigma =u_1\ldots  u_n\ldots  \in \Si^\om \mid  \fa i\geq 1 ~~ u_i\in V \}$
is the $\om$-power of $V$.
 The concatenation product is extended to the product of a 
finite word $u$ and an $\om$-word $v$: 
the infinite word $u.v$ is then the $\om$-word such that:
 $(u.v)(k)=u(k)$  if $k\leq |u|$ , and  $(u.v)(k)=v(k-|u|)$  if $k>|u|$.

\hs If $A$ is a subset of $B$ we shall denote $A^-=B-A$ the complement of $A$ (in $B$). 

\hs We assume the reader to be familiar with the theory of formal languages and of 
\orl s. We recall that \orl s form the class of \ol s accepted 
by finite automata with a   B\"uchi acceptance condition and this class is the omega Kleene 
closure of the class of regular finitary languages. 

\hs 
 We are going now to introduce  the notion of infinitary rational relation 
$R \subseteq \Si_1^\om \times \Si_2^\om$ via acceptance by 
2-tape B\"uchi automata. 

\begin{Deff}
\noi A  2-tape B\"uchi automaton 
 is a 7-tuple $\mathcal{T}=(K, \Si_1, \Si_2,   \Delta, q_0, F)$, where 
$K$ is a finite set of states, $\Si_1$, $\Si_2$,   are finite  alphabets, 
$\Delta$ is a finite subset of 
$K \times \Si_1^\star  \times \Si_2^\star  \times K$ 
called the set of transitions, $q_0$ is the initial state,  and $F \subseteq K$ is the set of 
accepting states. 
\nl A computation $\mathcal{C}$ of the  
2-tape B\"uchi automaton $\mathcal{T}$ over the pair 
$(u, v) \in \Si_1^\om \times \Si_2^\om$ 
 is an infinite sequence of transitions 
$$(q_0, u_1, v_1,  q_1), (q_1, u_2, v_2, q_2), \ldots 
(q_{i-1}, u_{i}, v_{i},   q_{i}), 
(q_i, u_{i+1}, v_{i+1}, q_{i+1}), \ldots $$
\noi such that:   $u=u_1.u_2.u_3 \ldots$ ~~and ~~  
 $v=v_1.v_2.v_3 \ldots$.
\nl The computation is said to be successful iff there exists an accepting state $q_f \in F$ 
and infinitely many integers $i\geq 0$ such that $q_i=q_f$. 
\nl The infinitary rational relation 
$R(\mathcal{T})\subseteq \Si_1^\om \times \Si_2^\om$ 
accepted by the 2-tape B\"uchi automaton $\mathcal{T}$ 
is the set of pairs $(u, v) \in \Si_1^\om \times \Si_2^\om$ 
such that there is some successful computation $\mathcal{C}$ of $\mathcal{T}$ over $(u, v)$. 
\nl The set of infinitary rational relations accepted by 
2-tape B\"uchi automata 
will be denoted $RAT_2$. 

\end{Deff} 

\noi As noticed in the introduction an infinitary rational relation 
$R \subseteq \Si_1^\om \times \Si_2^\om$ may be considered as an  $\om$-language 
over the product alphabet $\Si_1 \times \Si_2$.
 We shall use this fact to investigate the topological complexity 
of infinitary rational relations.

\section{Borel sets}

\noi We assume the reader to be familiar with basic notions of topology which
may be found in  \cite{mos,kec,lt,sta,pp}.

\hs For a  finite alphabet $X$ 
 we shall consider $X^\om$ as a topological space with the Cantor topology.
 The open sets of $X^\om$ are the sets in the form $W.X^\om$, where $W\subseteq X^\star$.
A set $L\subseteq X^\om$ is a closed set iff its complement $X^\om - L$ is an open set.
\nl  Define now the next classes of the  Hierarchy of Borel sets of finite ranks:

\begin{Deff}
The classes ${\bf \Si_n^0}$ and ${\bf \Pi_n^0 }$ of the Borel Hierarchy
 on the topological space $X^\om$  are defined as follows:
\nl ${\bf \Si^0_1 }$ is the class of open sets of $X^\om$.
\nl ${\bf \Pi^0_1 }$ is the class of closed sets of $X^\om$.
\nl And for any integer $n\geq 1$:
\nl ${\bf \Si^0_{n+1} }$   is the class of countable unions 
of ${\bf \Pi^0_n }$-subsets of  $X^\om$.
\nl ${\bf \Pi^0_{n+1} }$ is the class of countable intersections of 
${\bf \Si^0_n}$-subsets of $X^\om$.

\end{Deff}

\noi  The Borel Hierarchy is also defined for transfinite levels, but we shall not 
need them in the present study. 
There are also some subsets of $X^\om$ which are not Borel.  In particular 
the class of Borel subsets of $X^\om$ is strictly included into 
the class  ${\bf \Si^1_1}$ of analytic sets which are 
obtained by projection of Borel sets, 
see for example \cite{sta,lt,pp,kec} for more details.

\hs Recall also the notion of completeness with regard to reduction by continuous functions. 
For an integer $n\geq 1$, a set $F\subseteq X^\om$ is said to be 
a ${\bf \Si^0_n}$  (respectively,  ${\bf \Pi^0_n}$, ${\bf \Si^1_1}$)-complete set 
iff for any set $E\subseteq Y^\om$  (with $Y$ a finite alphabet): 
 $E\in {\bf \Si^0_n}$ (respectively,  $E\in {\bf \Pi^0_n}$,  $E\in {\bf \Si^1_1}$) 
iff there exists a continuous function $f: Y^\om \ra X^\om$ such that $E = f^{-1}(F)$.  
\nl  A ${\bf \Si^0_n}$
 (respectively,  ${\bf \Pi^0_n}$, ${\bf \Si^1_1}$)-complete set is a ${\bf \Si^0_n}$
 (respectively,  ${\bf \Pi^0_n}$, ${\bf \Si^1_1}$)-set which is 
in some sense a set of the highest 
topological complexity among the ${\bf \Si^0_n}$
 (respectively,  ${\bf \Pi^0_n}$, ${\bf \Si^1_1}$)-sets. 
 ${\bf \Si^0_n}$
 (respectively,  ${\bf \Pi^0_n}$)-complete sets, with $n$ an integer $\geq 1$, 
 are thoroughly characterized in \cite{stac}.  

  \begin{Exa}\label{exa} Let $\Si=\{0, 1\}$ and $\mathcal{A}=(0^\star.1)^\om \subseteq \Sio$. 
$\mathcal{A}$ is the set of 
$\om$-words over the alphabet $\Si$ with infinitely many occurrences of the letter $1$. 
It is well known that $\mathcal{A}$ is a 
${\bf \Pi^0_2 }$-complete set and its complement $\mathcal{A}^-$ 
is a ${\bf \Si^0_2 }$-complete set: it is the set of $\om$-words over $\{0, 1\}$ having 
only a finite number of occurrences of letter $1$. 
\end{Exa}

\section{${\bf \Pi^0_3}$-complete infinitary  rational relations}\label{pi}

\noi We had got in \cite{relratbor} some 
${\bf \Pi^0_3}$-complete infinitary  rational relations. 
We used a coding of $\om^2$-words over a finite alphabet $\Si$ by 
pairs of $\om$-words over $\Si \cup \{A\}$ where $A$ is an additionnal letter 
not in $\Si$.  
\nl We shall modify the previous proof (only sketched in \cite{relratbor}) 
by coding an $\om^2$-word over a finite alphabet $\Si$ 
by a {\bf single} $\om$-word over $\Si \cup \{A\}$. This way we can get some 
${\bf \Pi^0_3}$-complete infinitary  rational relation having some extra property. 

\begin{The}\label{Pi30} Let $\Ga=\{0, 1, A\}$ be an alphabet having three letters, and  
$\alpha$ be the $\om$-word over  the alphabet 
$\Ga$ which is defined by:

$$\alpha = A.0.A.0^2.A.0^3.A.0^4.A.0^5.A \ldots A.0^n.A.0^{n+1}.A \ldots$$
\noi Then there exists an 
infinitary  rational relation 
$R \subseteq \Ga^\om \times \Ga^\om$
such that:  
\nl $R_\alpha = \{\sigma \in \Ga^\om \mid (\sigma, \alpha) \in R \}  
\mbox{ is a } {\bf \Pi^0_3}-\mbox{complete subset of }\Ga^\om$,  and
 for all $u\in \Ga^\om-\{\alpha\}$~~
$R_u = \{\sigma \in \Ga^\om \mid (\sigma, u) \in R \} = \Ga^\om$. 
 Moreover $R$ is a  ${\bf \Pi^0_3}$-complete subset of $\Ga^\om \times \Ga^\om$. 

\end{The}

 \proo   We shall use a well known example of ${\bf \Pi^0_3}$-complete set which is a 
subset of the topological space $\Si^{\om^2}$. 

\hs The set  $\Si^{\om^2}$ is the set of $\om^2$-words over the finite alphabet $\Si$. 
It may also be viewed as the set of (infinite) $(\om \times \om)$-matrices 
whose coefficients are letters 
of $\Si$. If $x \in \Si^{\om^2}$ we shall write $x = (x(m, n))_{m\geq 1, n\geq 1}$.  
The infinite word $x(m, 1)x(m, 2)\ldots x(m, n)\ldots$ will be called 
the $m^{th}$ column of the $\om^2$-word $x$ and the infinite word 
$x(1, n)x(2, n)\ldots x(m, n)\ldots$ will be called 
the $n^{th}$ row of the $\om^2$-word $x$.  Thus an element of   $\Si^{\om^2}$ 
is completely determined by the (infinite) set of its columns or of its rows. 
\nl The set $\Si^{\om^2}$ is usually equipped with the product topology  of the dicrete 
topology on $\Si$ (for which every subset of $\Si$ is an open set), see \cite{kec} \cite{pp}. 
This topology may be defined 
by the following distance $d$. Let $x$ and $y$ be two  $\om^2$-words in $\Si^{\om^2}$ 
such that $x\neq y$, then  
$$ d(x, y)=\frac{1}{2^n} ~~~~~~~\mbox{   where  }$$
$$n=min\{p\geq 1 \mid \exists (i, j) ~x(i, j)\neq y(i, j) \mbox{ and } i+j=p\}$$

\hs Then the topological space $\Si^{\om^2}$ is homeomorphic to the above defined 
topological space  $\Si^{\om}$.  
The Borel hierarchy and the projective hierarchy on  $\Si^{\om^2}$ are defined from open 
sets in the same manner as in the case of the topological space $\Si^\om$. The notion 
of ${\bf \Si^0_n}$ (respectively ${\bf \Pi^0_n}$)-complete sets are also 
defined in a similar way.  

\hs Let now 
$$P = \{ x\in \{0, 1\}^{\om^2} \mid  \fa m \exists^{<\infty} n ~x(m, n)=1\}$$ 
\noi where $\exists^{<\infty}$ means ``there exist only finitely many", 
\nl $P$ is the set of $\om^2$-words 
having  all their columns in the ${\bf \Si^0_2 }$-complete subset $\mathcal{A}^-$ 
of $\{0, 1\}^{\om}$ where $\mathcal{A}$ is the ${\bf \Pi^0_2 }$-complete \orl~
given in Example \ref{exa}. 

\hs Recall the following classical result, \cite[p. 179]{kec}: 

\begin{Lem}
The set $P$ 
is a ${\bf \Pi^0_3}$-complete subset of $\{0, 1\}^{\om^2}$. 
\end{Lem} 

\proo Let 
$\mathcal{B}_m = \{ x\in \Si^{\om^2} \mid  x(m,1)x(m,2)\ldots x(m,n)\ldots \in \mathcal{A}^-\}$ 
be the  set of $\om^2$-words over $\Si=\{0, 1\}$ 
having their $m^{th}$ column in the ${\bf \Si^0_2 }$-complete set $\mathcal{A}^-$. 
 In order to prove  that,  for every integer $m\geq 1$,  
the set $\mathcal{B}_m$ is a  ${\bf \Si^0_2}$-subset of $\Si^{\om^2}$,  
 consider the function 
$i_m: \Si^{\om^2} \ra \Sio$ defined by $i_m(x)=x(m,1)x(m,2)\ldots x(m,n)\ldots $ for 
every $x\in \Si^{\om^2}$. The function  $i_m$ is continuous  and 
$i_m^{-1}(\mathcal{A}^-)=\mathcal{B}_m$ holds. Therefore  $\mathcal{B}_m$ is a 
 ${\bf \Si^0_2}$-subset of $\Si^{\om^2}$ because the class ${\bf \Si^0_2}$ is closed 
under inverse images  by  continuous functions. 

\hs   Thus the set 
$$P=\bigcap_{m\geq 1}\mathcal{B}_m$$ 
\noi of $\om^2$-words over $\Si$ 
having all their columns in $\mathcal{A}^-$ is a countable intersection 
 of ${\bf \Si^0_2}$-sets so it is a 
${\bf \Pi^0_3}$-set. 

 \hs  It remains to show that $P$ is ${\bf \Pi^0_3}${\bf-complete}. Let then 
$L$ be a ${\bf \Pi^0_3}$-subset of $\Sio$. We know that $L=\cap_{i\in \mathbb{N}^\star}A_i$ 
for some ${\bf \Si^0_2}$-subsets 
$A_i$, $i\geq 1$, 
of $\Sio$. 
But  $\mathcal{A}^-$ is ${\bf \Si^0_2 }$-complete therefore, for each integer $i \geq 1$, 
there is some continuous function $f_i : \Sio \ra \Sio$ such that 
$f_i^{-1}(\mathcal{A}^-)=A_i$.

 \hs  Let  now $f$ be the  function  from $\Sio$ into $\Si^{\om^2}$ which is defined by 
$f(x)(m,n)=f_m(x)(n)$. The function $f$ is continuous because each function $f_i$ 
is continuous.

\hs   For $x\in \Sio$  $f(x)\in P$ iff the $\om^2$-word $f(x)$ has all its 
 columns in the \ol~ $\mathcal{A}^-$, i.e. iff for all integers $m \geq 1$~~ 
$$f_m(x)=f_m(x)(1)f_m(x)(2)\ldots f_m(x)(n) \ldots   \in \mathcal{A}^-$$
iff  $\fa m \geq 1 ~~ x\in A_m$. Thus  $f(x)\in  P$ iff 
$x\in L=\cap_{m\geq 1}A_m$ so $L=f^{-1}(P)$. 

\hs   We have then proved that 
all ${\bf \Pi^0_3 }$-subsets of $\Sio$ are inverse images by continuous functions of 
the  ${\bf \Pi^0_3 }$-set $P$ therefore $P$ is 
a ${\bf \Pi^0_3 }$-complete set.         \ep 

\hs In order to use this example we shall firstly define a coding of $\om^2$-words 
over $\Si$ by $\om$-words over the alphabet  $(\Si\cup\{A\})$ where  
$A$ is a new letter not in $\Si$.

\hs Let us call, 
for $x\in \Si^{\om^2}$ and $p$ an integer $\geq 1$:  
$$T^x_{p+1}=\{x(p,1), x(p-1, 2), \ldots , x(2, p-1), x(1,p)\}$$ 
the set of elements $x(m, n)$ with $m+n=p+1$ and 
$$U^x_{p+1}=x(p,1).x(p-1, 2) \ldots x(2, p-1).x(1,p)$$ 
the sequence formed by the concatenation of elements $x(m, n)$ of $T^x_{p+1}$ for 
increasing values of $n$.

\hs We shall code an $\om^2$-word $x\in \Si^{\om^2}$ by the $\om$-word $h(x)$ defined by
$$h(x)=A.U^x_2.A.U^x_3.A.U^x_4.A.U^x_5.A.U^x_6.A \ldots A.U^x_{n}.A.U^x_{n+1}.A\ldots $$

\noi Let then $h$ be the mapping from $\Si^{\om^2}$  into
 $(\Si\cup\{A\})^\om$ 
such that,  for every $\om^2$-word $x$ over the  alphabet $\Si$,  
$h(x)$ is the code  of the $\om^2$-word $x$ as defined above.
It is easy to see, from the definition of $h$ and of the  order of the enumeration 
of letters $x(m, n)$ in $h(x)$ (they are enumerated for  increasing values of $m+n$), 
that $h$ is a continuous function from $\Si^{\om^2}$  into 
$(\Si\cup\{A\})^\om$. 

\hs  Remark that the above coding of $\om^2$-words resembles the use of the Cantor pairing 
function as it was used to construct the complete sets $P_i$ and $S_i$ in \cite{sw} 
(see also \cite{stac} or \cite[section 3.4]{sta}). 

\begin{Lem}\label{lem}
Let $\Si$ be a finite alphabet. If $L \subseteq \Si^{\om^2}$ is  
${\bf \Pi^0_3}$-complete then 
$$h(L) \cup h(\Si^{\om^2})^-$$
\noi  is a ${\bf \Pi^0_3}$-complete subset of $(\Si \cup\{A\})^\om$.
\end{Lem}

\proo  The topological space $\Si^{\om^2}$  is compact 
thus its image by the continuous function 
$h$ is also a compact subset of the topological space 
$(\Si \cup\{A\})^\om$. 
The set  $h(\Si^{\om^2})$ is compact hence  it is a closed subset of 
$(\Si \cup\{A\})^\om$ and  its complement 
$$(h(\Si^{\om^2}))^- =  (\Si \cup\{A\})^\om - h(\Si^{\om^2})$$
\noi  is an open (i.e. a ${\bf \Si^0_1}$) subset of $(\Si \cup\{A\})^\om$.

\hs  On the other hand the function $h$ is also injective 
thus it is a bijection from $\Si^{\om^2}$  onto 
$h(\Si^{\om^2})$. But a continuous bijection between two compact sets is an homeomorphism
therefore $h$ induces an homeomorphism between  $\Si^{\om^2}$  and  $h(\Si^{\om^2})$. 
By hypothesis $L$ is a ${\bf \Pi^0_3}$-subset of $\Si^{\om^2}$  thus 
$h(L)$ is a 
${\bf \Pi^0_3}$-subset of  $h(\Si^{\om^2})$ (where Borel sets of the topological 
space $h(\Si^{\om^2})$ are defined from open sets as in the cases of the topological 
spaces $\Sio$ or $\Si^{\om^2}$). 

\hs  The topological space $h(\Si^{\om^2})$ is a 
topological subspace of $ (\Si \cup\{A\})^\om$ and its 
topology  is induced by the topology on $(\Si \cup\{A\})^\om$: open sets 
of $h(\Si^{\om^2})$ are traces on $h(\Si^{\om^2})$ of open sets of 
$(\Si \cup\{A\})^\om$ and the same result holds for closed sets. Then 
one can easily show by induction that for every integer  $n \geq 1$,  
${\bf \Pi^0_n  }$-subsets 
(resp.  ${\bf \Si^0_n  }$-subsets) of  $h(\Si^{\om^2})$ are traces on $h(\Si^{\om^2})$ of 
${\bf \Pi^0_n  }$-subsets 
(resp.  ${\bf \Si^0_n  }$-subsets) of  $(\Si \cup\{A\})^\om$, i.e. are 
intersections with $h(\Si^{\om^2})$ of 
${\bf \Pi^0_n  }$-subsets 
(resp.  ${\bf \Si^0_n  }$-subsets) of  $(\Si \cup\{A\})^\om$. 

\hs  But  $h(L)$ is a ${\bf \Pi^0_3}$-subset of  $h(\Si^{\om^2})$ hence there exists 
a  ${\bf \Pi^0_3}$-subset $T$ of  $(\Si \cup\{A\})^\om$ such that 
$h(L)=T \cap h(\Si^{\om^2})$. But $h(\Si^{\om^2})$ is a closed 
i.e. ${\bf \Pi^0_1}$-subset (hence also a ${\bf \Pi^0_3}$-subset) 
of $(\Si \cup\{A\})^\om$  and the class of ${\bf \Pi^0_3}$-subsets of 
$(\Si \cup\{A\})^\om$  is closed under finite intersection thus 
$h(L)$ is a ${\bf \Pi^0_3}$-subset of $(\Si \cup\{A\})^\om$.  
 
\hs  Now  $h(L) \cup  (h(\Si^{\om^2}))^-$ is the union of a ${\bf \Pi^0_3}$-subset 
and of a ${\bf \Si^0_1}$-subset of $(\Si \cup\{A\})^\om$ 
therefore it is a ${\bf \Pi^0_3}$-subset of $(\Si \cup\{A\})^\om$ 
because the class of 
${\bf \Pi^0_3}$-subsets of $(\Si \cup\{A\})^\om$ 
 is closed under finite union.  

\hs  In order to prove that $h(L) \cup  (h(\Si^{\om^2}))^-$ 
 is ${\bf \Pi^0_3}${\bf-complete} it suffices to remark 
that
$$L=h^{-1}[ h(L) \cup  (h(\Si^{\om^2}))^- ]$$
\noi This implies that $h(L) \cup  (h(\Si^{\om^2}))^- $ is ${\bf \Pi^0_3}$-complete 
because  $L$ is assumed to be ${\bf \Pi^0_3}$-complete.    \ep 

\begin{Lem}\label{lem2}
Let  $P = \{ x\in \{0, 1\}^{\om^2} \mid  \fa m \exists^{<\infty} n ~x(m, n)=1\}$
 and  $\Si=\{0, 1\}$. 
Then  
$$\mathcal{P}=  h(P) \cup (h(\Si^{\om^2}))^- $$ 
\noi  is a ${\bf \Pi^0_3}$-complete subset of $(\Si \cup\{A\})^\om$.
\end{Lem}

\proo It follows directly from the two preceding Lemmas. \ep 

\hs Let now $\Si=\{0, 1\}$ and let $\alpha$ be the $\om$-word over the alphabet 
$\Si \cup\{A\}$ which is defined by:

$$\alpha = A.0.A.0^2.A.0^3.A.0^4.A.0^5.A \ldots A.0^n.A.0^{n+1}.A \ldots$$

\noi We can now state the following Lemma.

\begin{Lem}\label{R1}
Let $\Si=\{0, 1\}$ and $\alpha$ be the $\om$-word over  $\Si \cup\{A\}$ defined as above. 
Then there exists  an infinitary rational relation 
$R_1 \subseteq (\Si\cup\{A\})^\om \times (\Si \cup\{A\})^\om$ such that:
$$\fa x\in \Si^{\om^2}~~~ (x\in P) \mbox{  iff } ( (h(x), \alpha) \in R_1 )$$ 
\end{Lem}

\proo  We define now the relation $R_1$.
A pair 
$y=(y_1, y_2)$ of $\om$-words over the alphabet $\Si\cup\{A\}$ is in $R_1$ 
if and only if it is in the form

\hs $y_1 = U_k.u_1.v_1.A.u_2.v_2.A.u_3.v_3.A  
\ldots A.u_{n}.v_{n}.A. \ldots$
\nl $y_2 = V_k.w_1.z_1.A.w_2.z_2.A.w_3.z_3.A \ldots 
 A.w_{n}.z_{n}.A \ldots$

\hs where $k$ is an integer $\geq 1$, $U_k, V_k  \in (\Sis.A)^k$,  and,  
for all integers $i\geq 1$, 
$$v_i, w_i, z_i \in 0^\star \mbox{ and }  u_i \in \Sis  \mbox{  and   }  $$ 
$$|w_i|=|v_i|  ~~~~ \mbox{ and  }  ~~~~[~~ |u_{i+1}|=|z_i|+1 \mbox{ or } |u_{i+1}|=|z_i| ~~]$$ 
and there exist infinitely many integers $i$ such that $|u_{i+1}|=|z_i|$.

\hs We  prove first that the  relation  $R_1$ satisfies:   
$$\fa x\in \Si^{\om^2}~~~ (x\in P) \mbox{  iff } ( (h(x), \alpha) \in R_1 )$$ 

\noi Assume that for some $x\in \Si^{\om^2}$~~~ $(h(x), \alpha) \in R_1$. Then 
$(h(x), \alpha)$ may be written in the above form $(y_1,y_2)$ with 
\hs $y_1 = U_k.u_1.v_1.A.u_2.v_2.A.u_3.v_3.A  
\ldots A.u_{n}.v_{n}.A. \ldots$
\nl $y_2 = V_k.w_1.z_1.A.w_2.z_2.A.w_3.z_3.A \ldots 
 A.w_{n}.z_{n}.A \ldots$

\hs $y_1=h(x)$ implies that for all integers $n\geq 1$ ~~ $U_{k+n}^x=u_n.v_n$ thus 
$|u_n.v_n|=k+n-1$. 
\nl $y_2=\alpha$  implies that for all integers $n\geq 1$ ~~ $w_n.z_n=0^{k+n-1}$ thus 
$|w_n.z_n|=k+n-1$. 

\hs So  $|u_n.v_n|=|w_n.z_n|$ but by hypothesis $|w_n|=|v_n|$ therefore $|u_n|=|z_n|$. 

\hs Moreover $|u_{n+1}|=|z_n|+1$ or $|u_{n+1}|=|z_n|$. 
\nl If $|u_{n+1}|=|z_n|+1$ then $|u_{n+1}|=|u_n|+1$ and $|v_{n+1}|=|v_n|$ because 
$|u_{n+1}|+|v_{n+1}|=|u_{n}|+|v_{n}|+1$. 
\nl If $|u_{n+1}|=|z_n|$ then $|u_{n+1}|=|u_n|$ and $|v_{n+1}|=|v_n|+1$ because 
$|u_{n+1}|+|v_{n+1}|=|u_{n}|+|v_{n}|+1$. 

\hs This proves that the sequence $(|v_n|)_{n\geq1}$ is increasing because 
for all integers $n\geq 1$ ~~ $|v_{n+1}|=|v_n|$ or $|v_{n+1}|=|v_n|+1$. Moreover 
by definition of $R_1$ we know that there exist infinitely many integers $n\geq 1$ such that 
$|u_{n+1}|=|z_n|$ hence also $|v_{n+1}|=|v_n|+1$.  
 Thus 
$$\lim_{n\ra + \infty} |v_n| = +\infty $$

\noi Let now $K$ be an integer $\geq 1$ and let us prove that the $K$ first columns of the 
$\om^2$-word $x$ have only finitely many occurrences of the letter $1$. 

\hs $\lim_{n\ra + \infty} |v_n| = +\infty $ thus there exists an integer $N\geq 1$ such that 
$\fa n\geq N$ ~~ $|v_n|\geq K$. 

\hs Consider now, for  $ n\geq N$,  
$$U_{k+n}^x=u_n.v_n=x(k+n-1,1).x(k+n-2, 2) \ldots x(2, k+n-2).x(1,k+n-1)$$ 

\noi We know that $v_n \in 0^\star$ thus 
$$x(|v_n|,k+n-|v_n|)=x(|v_n|-1, k+n+1-|v_n|)= \ldots =x(2, k+n-2)=x(1,k+n-1)=0$$ 
\noi and in particular 
$$x(K,k+n-K)=x(K-1, k+n+1-K)= \ldots =x(2, k+n-2)=x(1,k+n-1)=0$$ 
\noi because $|v_n|\geq K$. 

\hs These equalities hold for all integers $ n\geq N$ and this proves that the $K$ first 
columns of the $\om^2$-word $x$ have only finitely many occurrences of the letter $1$. 
\nl But this is true for all integers $K \geq 1$ so {\bf all columns } of $x$ have  
a finite number of occurrences of the 
letter $1$ and $x\in P$. 

\hs Conversely it is easy to see that for each  $x\in P$ the pair 
 $(h(x), \alpha)$ may be written 
in the above form $(y_1, y_2) \in R_1$.  

\hs It remains only to prove that the above defined relation $R_1$ is an infinitary rational 
relation.  It is easy to see that the following 
$2$-tape  B\"uchi automaton $\mathcal{T}$ accepts 
the infinitary rational relation $R_1$. 

\hs  $\mathcal{T}=(K,  \Ga, \Ga, \delta, q_0, F)$, 
where $K=\{q_0, q_1, q_2, q_3,  q_4, q_5\}$
is a finite set of states, $\Ga=\Si\cup\{A\}=\{0, 1, A\}$, 
with $\Si=\{0, 1\}$, 
$q_0$ is the initial state,  and $F = \{q_4\}$ is the set of 
final states. 
Moreover $\delta \subseteq K \times  \Ga^\star \times \Ga^\star \times K$ is the finite  
 set of transitions, containing  the following transitions:

\hs  $(q_0, a, \lambda,  q_0)$, for all $a\in\Si$, 
\nl  $(q_0, \lambda, a,  q_0)$, for all $a\in\Si$,
\nl  $(q_0, A, A,  q_0)$, 
\nl  $(q_0, A, A,  q_1)$, 
\nl  $(q_1, a, \lambda,  q_1)$, for all $a\in\Si$, 
\nl  $(q_1, \lambda, \lambda,  q_2)$,
\nl  $(q_2, a, 0,  q_2)$, for all $a\in\Si$, 
\nl  $(q_2, A, \lambda, q_3)$, 
\nl  $(q_3, a, 0,  q_3)$, for all $a\in\Si$,
\nl  $(q_3,\lambda, \lambda, q_4)$
\nl  $(q_3, a, \lambda,  q_5)$, for all $a\in\Si$,  
\nl  $(q_4, \lambda, A,  q_2)$, 
\nl  $(q_5, \lambda, A,  q_2)$.

\begin{Rem} Using classical constructions from automata theory, we could  
have avoided  the set of transitions to contain 
some transitions in the form  
$(q_i, \lambda, \lambda,  q_j)$, like $(q_1, \lambda, \lambda,  q_2)$ or 
$(q_3,\lambda, \lambda, q_4)$. 
\end{Rem}

\begin{Lem}\label{complement} The set 
$$R_2 = (\Si\cup\{A\})^\om\times (\Si \cup\{A\})^\om - ( h(\Si^{\om^2}) \times \{\alpha\} )$$
\noi is an infinitary rational relation.  
\end{Lem}

\proo By definition of the mapping $h$, we know that a pair of $\om$-words 
over   the alphabet  $(\Si\cup\{A\})$ is in $h(\Si^{\om^2}) \times \{\alpha\}$ iff 
 it is  in the form $(\sigma_1, \sigma_2)$, 
where
\hs  $\sigma_1 = A.u_1.A.u_2.A.u_3.A.u_4.A \ldots .A.u_n.A.u_{n+1}.A \ldots$
\nl  $\sigma_2 = \alpha = A.0.A.0^2.A.0^3.A.0^4.A \ldots A.0^n.A.0^{n+1}.A \ldots$

\hs where for all integers $i\geq 1$, $u_i \in \Sis$ and $|u_i|=i$. 

\hs So it is easy to see that 
$(\Si\cup\{A\})^\om\times (\Si \cup\{A\})^\om - ( h(\Si^{\om^2}) \times \{\alpha\} )$
is the union of the sets $\mathcal{C}_j$ where:

\begin{itemize} 
\ite $\mathcal{C}_1 = \{ (\sigma_1, \sigma_2) ~/~ \sigma_1 , \sigma_2 \in (\Si\cup\{A\})^\om 
\mbox{ and } ( \sigma_1 \in \mathcal{B} \mbox{ or } \sigma_2 \in \mathcal{B} ) \}$
\nl where $\mathcal{B}$ is the set of $\om$-words over $(\Si\cup\{A\})$ having only 
a finite number of letters $A$. 

\ite $\mathcal{C}_2$ is formed by pairs  $(\sigma_1, \sigma_2)$ where 
\nl $\sigma_1$ or $\sigma_2$ has not any initial segment in $A.\Si.A.\Si^2.A$.

\ite $\mathcal{C}_3$ is formed by pairs  $(\sigma_1, \sigma_2)$ where 
\nl $\sigma_2 \notin \{0, A\}^\om$. 

\ite $\mathcal{C}_4$ is formed by pairs  $(\sigma_1, \sigma_2)$ where 
\nl $\sigma_1 = A.w_1.A.w_2.A.w_3.A.w_4.A \ldots A.w_n.A.u.A.z_1 $
\nl $\sigma_2 = A.w'_1.A.w'_2.A.w'_3.A.w'_4.A \ldots A.w'_n.A.v.A.z_2 $

\hs where $n$ is an integer $\geq 1$,  for all $i \leq n$~  $w_i, w'_i \in \Sis$, 
$z_1, z_2 \in (\Si\cup\{A\})^\om$ and 
$$u, v \in \Sis \mbox{  and }  |v|\neq |u|$$

\ite $\mathcal{C}_5$ is formed by pairs  $(\sigma_1, \sigma_2)$ where 
\nl $\sigma_1 = A.w_1.A.w_2.A.w_3.A.w_4 \ldots A.w_n.A.w_{n+1}.A.v.A.z_1 $
\nl $\sigma_2 = A.w'_1.A.w'_2.A.w'_3.A.w'_4 \ldots A.w'_n.A.u.A.z_2 $

\hs where $n$ is an integer $\geq 1$,  for all $i \leq n$~  $w_i, w'_i \in \Sis$, 
$w_{n+1} \in \Sis$, 
$z_1, z_2 \in (\Si\cup\{A\})^\om$ and 
$$u, v \in \Sis \mbox{  and }  |v|\neq |u|+1$$

\end{itemize}

\noi Each set $\mathcal{C}_j$, $1\leq j\leq 5$, is easily seen to be an infinitary 
rational relation $\subseteq (\Si\cup\{A\})^\om \times (\Si\cup\{A\})^\om$ (the detailed 
proof is left to the reader).  The class $RAT_2$ is closed under finite union thus

$$R_2 = (\Si\cup\{A\})^\om\times (\Si \cup\{A\})^\om - ( h(\Si^{\om^2}) \times \{\alpha\} )
 = \bigcup_{1\leq j\leq 5} \mathcal{C}_j$$

\noi is an infinitary rational relation. \ep

\hs Return now to the proof of Theorem \ref{Pi30}.  Let 
$$R = R_1 \cup R_2 \subseteq \Gao \times \Gao$$
\noi The class $RAT_2$ is closed under finite union therefore $R$ is an 
infinitary rational relation. 

\hs Lemma \ref{R1} and the definition of $R_2$ imply  that 
$R_\alpha = \{\sigma \in \Ga^\om \mid (\sigma, \alpha) \in R \}$ is equal to the set 
$\mathcal{P}=  h(P) \cup (h(\Si^{\om^2}))^- $ 
which   is a ${\bf \Pi^0_3}$-complete subset of $(\Si \cup\{A\})^\om$ 
by Lemma \ref{lem2}. 

\hs Moreover, for all $u\in \Ga^\om-\{\alpha\}$, 
$R_u = \{\sigma \in \Ga^\om \mid (\sigma, u) \in R \} = \Ga^\om$ holds by 
definition of $R_2$.       

\hs In order to prove that $R$ is a ${\bf \Pi^0_3}$-set remark first that 
 $R$ may be written as the  union: 
 $$R = \mathcal{P} \times \{\alpha\} ~~ \bigcup  ~~\Gao \times (\Gao - \{\alpha\})$$
\noi We already know that  $\mathcal{P}$ is a ${\bf \Pi^0_3}$-complete subset of 
$(\Si\cup\{A\})^\om$. Then it is easy to show that $\mathcal{P} \times \{\alpha\}$ 
is also a ${\bf \Pi^0_3}$-subset of $(\Si\cup\{A\})^\om \times (\Si\cup\{A\})^\om$. 
 On the other side it is easy to see that  $\Gao \times (\Gao - \{\alpha\})$ is an open  
subset of $\Gao \times \Gao$.
Thus $R$ is a ${\bf \Pi^0_3}$-set because the Borel class ${\bf \Pi^0_3}$ 
is closed under finite union. 
 
\hs Moreover let  $g: \Si^{\om^2} \ra 
 (\Si\cup\{A\})^\om \times (\Si\cup\{A\})^\om $ be the function defined  by:

$$\fa x \in \Si^{\om^2}~~~~~~~ g(x) = (h(x) , \alpha)$$

\noi  It is easy to see that $g$ is continuous because $h$ is continuous. By construction 
it turns out that 
for all $\om^2$-words $x\in \Si^{\om^2}$ ~~~ $(x\in P)$ iff $(g(x)\in R)$. 
This means that $g^{-1}(R)= P$. This implies that $R$ is ${\bf \Pi^0_3}$-complete 
because $P$ is ${\bf \Pi^0_3}$-complete. \ep

\begin{Rem} The structure of the 
${\bf \Pi^0_3}$-complete infinitary rational relation $R$ we have just got is very different 
from the structure of a previous  example given in \cite{relratbor}. It can be described 
very simply by the sections $R_u$, $u\in \Gao$.  All sections but one are equal to 
$\Gao$, so they have the lowest topological complexity and exactly one section is a 
${\bf \Pi^0_3}$-complete subset of $\Gao$. 
\end{Rem}

\hs {\bf  Acknowledgements.} Thanks to Jean-Pierre Ressayre and  Pierre Simonnet
for useful discussions.

\end{document}